\documentclass[12pt]{article}
\usepackage{graphics}
\usepackage{amssymb}

\newcommand{\e}{\mathrm{e}}
\newcommand{\N}{\mathbb{N}}
\newcommand{\R}{\mathbb{R}}
\newcommand{\Z}{\mathbb{Z}}
\newcommand{\BB}{\mathcal{B}}
\newcommand{\CC}{\mathcal{C}}
\newcommand{\LL}{\mathcal{L}}
\newcommand{\OO}{\mathcal{O}}
\newcommand{\PP}{\mathcal{P}}
\newcommand{\UU}{\mathcal{U}}
\newcommand{\pd}{\partial}
\newtheorem{claim}{Claim}[section]
\newtheorem{theorem}[claim]{Theorem}

\newtheorem{lemma}[claim]{Lemma}
\newtheorem{remark}[claim]{Remark}
\newtheorem{remarks}[claim]{Remarks}
\newenvironment{proof}[1][Proof]{\textsl{#1:} }{\ \rule{0.4em}{0.7em}}
\begin{document}

\title{Absolute continuity in periodic thin tubes \\ and
strongly coupled leaky wires}
\author{Francois Bentosela$^{1,2}$, Pierre Duclos$^{1,3}$, 
and Pavel Exner$^{4,5}$}
\date{}
\maketitle

\begin{quote}
\emph{1$\;$ Centre de Physique Th\'eorique, C.N.R.S., Luminy Case
907, \\ \phantom{aa} F-13288 Marseille Cedex 9, \\
2$\;$ Universit\'{e} de la Mediterran\'{e}e
(Aix--Marseille II), F-13288 \\ \phantom{aa} Marseille, \\
3$\;$ PhyMat, Universit{\'e} de
Toulon et du Var, BP 132, F-83957 \\ \phantom{aa} La Garde Cedex, France; \\
4$\;$ Department of Theoretical Physics, Nuclear Physics Institute, \\
 \phantom{aa} Academy of Sciences, CZ-25068 \v{R}e\v{z} near Prague, \\
5$\;$ Doppler Institute, Czech Technical University, B\v{r}ehov{\'a}
7, \\ \phantom{aa} CZ-11519 Prague, Czechia; \\
\texttt{\phantom{aa} bento@cpt.univ-mrs.fr, duclos@univ-tln.fr, \\ 
\phantom{aa} exner@ujf.cas.cz}} \\[5mm]
{\small {\bf Abstract:}  Using a perturbative argument, we show
that in any finite region containing the lowest transverse eigenmode,
the spectrum of a periodically
curved smooth Dirichlet tube in two or three dimensions is
absolutely continuous provided the tube is sufficiently thin. In a
similar way we demonstrate absolute continuity at the bottom of
the spectrum for generalized Schr\"odinger operators with a
sufficiently strongly attractive $\delta$ interaction supported by
a periodic curve in $\mathbb{R}^d,\: d=2,3$.}
\end{quote}


\section{Introduction}

In models of periodically structured quantum systems, absolute
continuity of the spectrum is a crucial property. For usual
Schr\"odinger operators and many other PDE's with periodic
coefficients the problem is well understood -- see, e.g. \cite{Ku,
RS}. On the other hand, there are important classes of operators
which still pose open questions. An example is represented by
so-called \emph{quantum waveguides}, i.e. systems the Hamiltonian
of which is (a multiple of) the Laplacian in an infinitely long
tube-shaped region, usually with Dirichlet boundary conditions.

In the two-dimensional setting, where the region in question is a
periodically curved planar strip the absolute continuity has been
demonstrated recently in \cite{SV}. Unfortunately, the method used
in this work does not seem to generalize to other dimensions
including the physically interesting case of a periodic tube in
$\R^d,\, d=3$. This is why we present in this letter a simpler
result stating the absolute continuity at the bottom of the
spectrum for tubes which are thin enough. With physical
applications in mind we formulate it for $d=2,3$, but the argument
can be used in any dimension.

We also address an analogous problem concerning Schr\"odinger
operators in $L^2(\R^d)$, $d=2,3$, with an attractive $\delta$
interaction supported by a periodic curve; we prove absolute
continuity, again at the bottom of the spectrum, for a
sufficiently strong attraction. If $d=2$ the answer was known for
a family of curves periodic in two independent directions
\cite{BSS}, however, for a single infinite curve results have been
missing, to say nothing about the more singular case of dimension
$d=3$.

Our method is perturbative. We show that for a sufficiently thin
tube or strongly attractive $\delta$ interaction the Floquet
eigenvalues do not differ much from the Floquet eigenvalues of a suitable
comparison problem in one dimension which are known to be
nonconstant as functions of quasimomentum. A similar argument was
used recently for periodically perturbed magnetic channels
\cite{EJK}, and in a different context to demonstrate existence of
persistent currents in leaky quantum wire loops \cite{EY2}.


\section{Thin curved tubes}
\setcounter{equation}{0}

Let $\Gamma:\: \R\to\R^d,\, d=2,3$, be a $C^4$ smooth curve
without self-intersections which is periodic, i.e. there are $L>0$
and a nonzero vector $b\in\R^d$ such that
$$ 
  \Gamma(s+L) = b+\Gamma(s)\,, \quad\forall s\in\R\,.
$$ 
With an abuse of notation we will use the same symbol $\Gamma$ for
the map and for the image $\Gamma(\R)$. Let further
$\Omega:= \{x\in\R^d:\, \mathrm{dist}\,(x,\Gamma)<a\}$
for a fixed $a>0$. If the latter is small enough, such a tube has
no self-intersections and we can parametrize it by natural
curvilinear coordinates. More specifically, take a ball
$\BB_a:=\{u\in\R^{d-1},\,|u|<a\}$ and consider the maps $\phi:\:
\R\times\BB_a \to\R^d$ defined by
\begin{eqnarray*}
  \phi_\Gamma(s,u) &\!=\!& \Gamma(s) -un(s)\,, \\
  \phi_\Gamma(s,r,\vartheta) &\!=\!& \Gamma(s) -r[n(s)
  \cos(\vartheta\!-\!\beta(s)) + b(s)
  \sin(\vartheta\!-\!\beta(s))]
\end{eqnarray*}
for $d=2,3$, respectively, where $n(s), b(s)$ are the normal and
binormal vector to $\Gamma$ at the point $s$ and the function
$\beta$ in the three-dimensional case (rotation with respect to
the Frenet frame) will be specified later. We will always assume
that $a$ is so small that the above maps are diffeomorphisms.

The object of our interest is the Dirichlet Laplacian
$-\Delta_D^{\Omega}$ defined conventionally
\cite[Sec.~XIII.15]{RS} as the self-adjoint operator associated
with the quadratic form $\psi\mapsto \|\nabla\psi\|^2
_{L^2(\Omega)}$ with the domain $W^{2,1}_0(\Omega)$.

The curvilinear coordinates allow us to ``straighten'' the tube,
i.e. to pass to an operator on $L^2(\CC_a)$, where $\CC_a:=
\R\times\BB_a$ is a straight cylinder. The appropriate unitary
operator $U:\: L^2(\Omega)\to L^2(\CC_a)$ is at that defined by
$U\psi:= g^{1/4}\psi\circ\phi_\Gamma$, where $g^{1/2}$ is the
corresponding Jacobian, $g^{1/2}= 1+u\gamma(s)$ for $d=2$ and
$g^{1/2}= 1+r\gamma(s) \cos(\vartheta\!-\!\beta(s))$ for $d=3$,
with $\gamma$ being the curvature of $\Gamma$. The price we pay
for the simplification of the region is the more complicated form
of the operator $H_{a,\Gamma}:= U(-\Delta_D^{\Omega})U^{-1}$; a
straightforward calculation \cite{ES, DE} gives
$$ 
H_{a,\Gamma} = -\partial_s h^{-2} \partial_s -\Delta_D^{\BB_a} +
V\,,
$$ 
where $h:=g^{1/2}$ and $V$ is the effective potential induced by
the geometry,
$$ 
V = -{\gamma^2\over 4h^2} + {h_{ss}\over 2h^3} - {5h_s^2\over
4h^4}\,.
$$ 
For $d=3$ the transformation requires to choose $\beta(s) =
\int_{s_0}^s \tau(s)\,ds$, where $\tau$ is the torsion of
$\Gamma$; the latter appear also in the effective potential coming
from the derivatives of $h$ by Frenet-Serret formulae
\cite{DE}.

As the first step we denote $\CC^L_a:= [0,L)\times\BB_a$ and
perform the usual Floquet decomposition over the Brillouin zone
$\BB:=[-\pi/L, \pi/L)$.
\begin{lemma} \label{tubefloq}
There is a unitary $\UU:\: L^2(\CC_a) \to \int^\oplus_\BB
L^2(\CC^L_a)\, d\theta$ such that
$$ 
\UU H_{a,\Gamma} \UU^{-1} = \int^\oplus_\BB H_{a,\Gamma}(\theta)\,
d\theta\,,
$$ 
where the fibre operator satisfies periodic b.c. in $s$ acting as
$$ 
H_{a,\Gamma}(\theta) = (-i\partial_s+\theta)
h^{-2} (-i\partial_s+\theta) -\Delta_D^{\BB_a} + V\,.
$$ 
\end{lemma}
\begin{proof}  This is a classical result. We use here the
modification, sometimes ascribed to Skriganov, of
\cite[Thm~XIII.88]{RS} with the Floquet-Bloch transform given for
all $(\theta,s,u)\in\BB\times[0,L)\times\BB_a$ by the formula
$$
(\UU f)(\theta,s,u):=\sum_{n\in\Z}\sqrt{L\over2\pi}\:
\mathrm{e}^{-in\theta L-i\theta s}f(s+L n,u)\,.
$$
An exhaustive discussion for $d=2$ can be found in \cite{Yo}, the
argument for $d=3$ is analogous.
 \end{proof} \vspace{1em}

We need also the character of $\theta$-dependence of the fibre
operators.

\begin{lemma} \label{tubeanal}
$\: \{ H_{a,\Gamma}(\theta):\: \theta\in\BB\, \}$ is a type A
analytic family.
\end{lemma}
\begin{proof}  Splitting $H_{a,\Gamma}(\theta)$ in two pieces
$H_{a,\Gamma}(0)=-\partial_s h^{-2} \partial_s -\Delta_D^{\BB_a} +
V$ and $H_{a,\Gamma}(\theta)-H_{a,\Gamma}(0)=
-i\theta(\partial_sh^{-2}+h^{-2}\partial_s)+\theta^2h^{-2}$, one
sees that the first piece is self-adjoint on $\{f\in
W^{2,2}(\CC_a^L),\, f(L,u)=f(0,u),\,\partial_sf(L,u)
=\partial_sf(0,u)\}$ and that the second one is entire analytic
and relatively bounded perturbation of $H_{a,\Gamma}(0)$ with zero
relative bound.
\end{proof}
\vspace{1em}

Our main tool is the perturbation theory w.r.t. $a$. Since the
argument follows closely \cite{DE} we will just sketch it. First
we use transverse scaling to pass to a unitarily equivalent
operator on $L^2(\CC^L_a)$ given by
 $$ 
 \tilde H_a(\theta):= (-i\pd_s+\theta) h_a^{-2} (-i\pd_s+\theta)
 -a^{-2}\Delta_D^{\BB_1} + V_a\,,
 $$ 
where $h_a(s,u):= h(s,au)$ and $h_a(s,r,\vartheta):=
h(s,ar,\vartheta)$ for $d=2,3$, respectively, and similarly for
$V_a$. Let $\chi_j$ and $\kappa_j^2$ be the eigenfuctions and
eigenvalues of $-\Delta_D^{\BB_1}$. Using this transverse basis we
pass to the matrix representation,
 $$ 
 \tilde H_{a,jk}(\theta)= a^{-2} \kappa_j^2\delta_{jk}+T_{jk}\,,\quad
 T_{jk}:= (-i\pd_s+\theta) (h_a^{-2})_{jk}(-i\pd_s+\theta)
 +V_{a,jk}\,,
 $$ 
where $f_{jk}:= \int_{-1}^1 f(\cdot,au) \,\chi_j(u)\chi_k(u)\, du$
if $d=2$ and similarly for $d=3$. We need a reference operator
which will be chosen in the form
 $$ 
 \tilde H_a^0(\theta):=I\otimes( -a^{-2}\Delta_D^{\BB_1}) + S(\theta)\otimes I\,,
 \quad S(\theta):= (-i\pd_s+\theta)^2-\,\frac{1}{4}\,\gamma^2
 $$ 
with periodic b.c. in the variable $s$. The spectrum of
$S(\theta)$ is purely discrete; we denote its eigenvalues arranged
in the ascending order as $\lambda_n(\theta)$ with the index
$n\in\N$. Recall that the spectrum of $S(\theta)$ is simple with a
possible exception of the endpoints of the Brillouin zone and
$\theta=0$. Let $K$ be a compact subset of $\BB$ which does not
contain the points $0$ and $\pm\pi/L$. The eigenvalues
$\epsilon_{jn}^0(a,\theta):= a^{-2}\kappa_j^2+ \lambda_n(\theta)$
of $\tilde H_a^0(\theta)$ are isolated and of finite multiplicity;
with the above choice of $K$ they even become simple on $K$ for
any fixed $n$ and $j=1$ if $a$ is small enough. They depend on
$a$, of course, but one can perform the perturbation expansion
with respect to $W(\theta):=H_a(\theta)-H_a^0(\theta)$ around
such running values. Mimicking the argument of \cite{DE} we come
to the following conclusion.

\begin{lemma} \label{tubepert}
Let $n_0\in\N$ and $E_{n_0}:=((2n_0\!-\!1)\pi/L)^2$. There exists
a positive $a_{K, n_0}$ such that for all  $a\in (0,a_{K, n_0})$
and any $\theta\in K$, the spectrum of $\tilde H_a(\theta)$ below
$E_{n_0}$ consists of exactly $n_0$ simple eigenvalues
$\{\epsilon_{1,n}(a,\theta)\}_{1\le n\le n_0}$. Moreover, the
expansion
$$
\epsilon_{1,n}(a,\theta)=a^{-2}\kappa_1^2+\lambda_n(\theta)+\OO(a)
$$
holds for each $n=1,\dots,n_0$ uniformly in $K$.
\end{lemma}
Armed with these
prerequisites, we can now formulate and prove the main result of
this section.

\begin{theorem} \label{tubemain}
To any $E>0$ there is $a_E>0$ such that the spectrum of
$-\Delta_D^{\Omega}$ is absolutely continuous in the interval
$[0,a^{-2}\kappa_1^2+E]$ for all $a<a_E$.
\end{theorem}
\begin{proof}
By Lemma~\ref{tubefloq} we have to check that the eigenvalues of
$\tilde H_a(\theta)$ are nowhere constant as functions of the
quasimomentum $\theta$. Since they are real-analytic by
Lemma~\ref{tubeanal}, we have only to verify that each eigenvalue
branch acquires at least two different values in the set $K$.
Without loss of generality we can put $E=E_{n_0}$ for some
$n_0\in\N$. Then there is just $n_0$ eigenvalues in
$[0,a^{-2}\kappa_1^2+E]$, and consequently, there exists a $c_E>0$
such that $|\epsilon_{1,n}(a,\theta) -a^{-2}\kappa_1^2
+\lambda_n(\theta)|\le c_E a$ holds for each of them according to
Lemma~\ref{tubepert}. Since the functions $\lambda_n(\cdot)$ are
non-constant by \cite[Sec.~XIII.16]{RS}, the conclusion follows.
\end{proof} \hspace{1em}

\begin{remark} \label{small gamma}
{\rm A similar perturbative argument shows that the spectrum is
absolutely continuous at its bottom if $a$ is fixed and
$\|\gamma\|_\infty$ is sufficiently small.}
\end{remark}


\section{Leaky quantum wires}
\setcounter{equation}{0}

Now we are going to consider the analogous problem for another
class of operators. Let $\Gamma$ be the same $C^4$-smooth periodic
curve in $\R^d$ and consider the generalized Schr\"odinger
operator given by the formal expression
 $$ 
 H_{\alpha,\Gamma} = -\Delta-\alpha \delta(x-\Gamma)\,.
 $$ 
Its meaning is different for different dimensions. If $d=2$ we
regard it as the unique self-adjoint operator associated with the
quadratic form
 $$ 
 q_{\alpha,\Gamma}[\psi] = \|\nabla\psi\|^2 - \alpha \int_\R
 |\psi(\Gamma(s))|^2 ds\,, \quad \psi\in W^{2,1}(\R^2)\,,
 $$ 
which is closed and below bounded by \cite{BEKS}; we suppose that
the singular interaction is attractive, $\alpha>0$. The situation
is more complicated in the three-dimensional case when
$\mathrm{codim}\,\Gamma=2$. Following the construction given in
\cite{EK1} -- see also \cite{Po} for a more general background --
one starts from the family of curves $\phi_\Gamma
(\cdot,\rho,\vartheta_0)$, obtained by translating $\Gamma$
by $\rho(\cos(\vartheta_0),\sin(\vartheta_0))$,  which are used
to determine the generalized
boundary values of a function $f\psi\in W_{\mathrm{loc}}^{2,2}(
\R^3 \setminus \Gamma)\cap L^{2}(\R^3)$ as the following limits
 \begin{eqnarray*}
 L_0(\psi)(s) &\!:=\!& -\lim_{\rho \to 0}\: \frac{1}{\ln \rho }\,
 \psi(\phi_\Gamma(s,\rho,\vartheta_0))\,, \\
 L_1(\psi)(s) &\!:=\!& \lim_{\rho \to 0}\,
 \left[\, \psi(\phi_\Gamma(s,\rho,\vartheta_0))
 +L_0(\psi)(s)\ln \rho \,\right] \,.
 \end{eqnarray*}
We call $\Upsilon_\Gamma$ the family of those $\psi$ for which
these limits exist a.e. in $\R$, are independent of $\vartheta_{0}$,
and define a pair functions belonging to $L_{\mathrm{loc}}^2(\R)$.
The sought operator is then defined as
 \begin{eqnarray*}
 H_{\alpha,\Gamma}\psi(x) &\!:=\!& -\Delta\psi(x)\,,\quad x\in
 \R^3 \setminus\Gamma\,, \\
 D(H_{\alpha,\Gamma}) &\!:=\!& \{f\in \Upsilon_\Gamma :\: 2\pi
 \alpha L_0(\psi)(s)=L_1(\psi)(s)\,\}\,.
 \end{eqnarray*}
Being interested in strong coupling we suppose hereafter that the
coupling parameter $\alpha$ is negative though $H_{\alpha,\Gamma}$
is well defined for all $\alpha\in\R$; recall that the
two-dimensional $\delta$ interaction is always attractive
\cite{AGHH}.

In an important particular case when $\Gamma$ is a straight line
one uses separation of variables to show that the spectrum is
absolutely continuous and covers the interval $[\zeta(\alpha),
\infty)$, with the threshold given by the corresponding
$(d\!-\!1)$-dimensional $\delta$ interaction eigenvalue,
 $$ 
 \zeta(\alpha) := \left\lbrace \begin{array}{lcl} -{1\over
 4} \alpha^2 \quad & \dots & \quad d=2 \\ -4\e^{2(-2\pi\alpha
 +\psi(1))} \quad & \dots & \quad d=3 \end{array} \right.
 $$ 
where $-\psi(1)\approx 0.577$ is the Euler number. This also
illustrates that the strong coupling means $(-1)^d\alpha
\to\infty$ for $d=2,3$. Due to the injectivity and periodicity
assumptions the curve decomposes into a disjoint union of
translates of the period cell $\Gamma_\PP:= \Gamma \upharpoonright
[0,L)$. Since $H_{\alpha,\Gamma}$ now acts in the whole Euclidean
space, we need also a decomposition of the space $\R^d$ with the
period cell
 $$ 
 \PP := \left\lbrace \, \LL+tb:\: t\in[0,1)\, \right\rbrace\,,
 $$ 
where $\LL\subset\R^d$ is a affine space which is not colinear
with $b$. We denote by $b_\perp$ the component of $b$ in the
direction orthogonal to $\LL$; it follows that $b_\perp\ne0$. It
is important that the two decompositions are chosen in a
consistent way, $\Gamma_\PP=\PP\cap\Gamma$. We will assume in
addition that
\\ [.5em]
 (c) the restriction of $\Gamma_\PP$ to the interior of $\PP$ is
 connected. \\ [.5em]
It should be noted that the choice of a slab for $\PP$ is made
rather for convenience -- see Remarks~\ref{puzzle,crochet} below.
We start again from the Floquet decomposition with respect to the
Brillouin zone $\BB:=[-\pi|b_\perp|^{-1}, \pi|b_\perp|^{-1})$.

\begin{lemma} \label{leakyfloq}
There is a unitary $\UU:\: L^2(\R^d) \to \int^\oplus_\BB
L^2(\PP)\, d\theta$ such that
$$ 
\UU H_{\alpha,\Gamma} \UU^{-1} = \int^\oplus_\BB
H_{\alpha,\Gamma}(\theta)\, d\theta\,,
$$ 
where the fibre operator satisfies periodic b.c. in the direction
of $b$ acting as
$$ 
H_{a,\Gamma}(\theta) = (-i\nabla+\theta)^2
-\alpha\delta(x-\Gamma)\,,
$$ 
 and the interaction term in $L^2(\PP)$ is interpreted in the
 above described sense, the quadratic form if $d=2$ and boundary
 conditions if $d=3$.
\end{lemma}
\begin{proof}  See \cite{EY1} for $d=2$ and \cite{EK2} for $d=3$.
\end{proof} \vspace{1em}

It is easy to see that in distinction to the previous case the
essential spectrum is non-empty and equals
$\sigma_\mathrm{ess}(H_{a,\Gamma}(\theta))= [\theta^2,\infty)$; we
will be interested in the eigenvalues below its threshold. They
are again real-analytic functions:

\begin{lemma} \label{leakyanal}
$\: \{ H_{\alpha,\Gamma}(\theta):\: \theta\in\BB\, \}$ is a type A
analytic family.
\end{lemma}
\begin{proof} Similar to that of Lemma~\ref{tubeanal}.
\end{proof} \vspace{1em}

The role of the small tube width from the previous section is
played here by strong coupling. While the wave function may be
nonzero at large distances from $\Gamma$, it is localized in its
vicinity as $(-1)^d\alpha \to\infty$. Then one can choose a
tubular neighbourhood $\Omega$ of $\Gamma$ and  estimate
the operator in question from both sides by imposing the Dirichlet
and Neumann condition at $\pd\Omega$. The exterior part does not
contribute to the negative spectrum, while the part in $\Omega$ can be
treated as in the previous section, with the additional $\delta$
interaction on the tube axis and different boundary conditions for
the lower bound. The argument is thus more complicated, however,
it was done in \cite{EY1} and \cite{EK2} with the following
result.

\begin{lemma} \label{leakypert}
The number of isolated eigenvalues of $H_{\alpha,\Gamma}(\theta)$
exceeds any fixed $n\in\N$ as $(-1)^d\alpha \to\infty$. The $n$-th
eigenvalue behaves asymptotically as
 $$ 
 \epsilon_n(\alpha,\theta):= \zeta(\alpha)+ \lambda_n(\theta) +
 \left\lbrace \begin{array}{l} \OO(\alpha^{-1}\ln\alpha) \\
 \OO(\e^{\pi\alpha}) \end{array} \right\rbrace
 $$ 
in the strong coupling limit for $d=2,3$, respectively, uniformly
in $\theta$, i.e. the error terms is for a fixed $n$ bounded in
$\BB$. Here $\lambda_n(\theta)$ means again the $n$-th eigenvalue
of the operator $S(\theta)$ defined in the previous section.
\end{lemma}

The main result of this section then reads as follows.

\begin{theorem} \label{leakymain}
Under the stated assumptions, to any $E>0$ there exists an
$\alpha_E>0$ such that the spectrum of the operator
$H_{\alpha,\Gamma}$ is absolutely continuous in $(-\infty,
\zeta(\alpha)+E]$ as long as $(-1)^d\alpha >\alpha_E$.
\end{theorem}
\begin{proof}
The argument is analogous to the one used for Theorem~\ref{tubemain}.
\end{proof} \hspace{1em}

\begin{remarks} \label{puzzle,crochet}
{\rm (i) It is not always possible to choose $\PP$ is the
Cartesian-product form, as we did above, which would satisfy the
assumption (c). Counterexamples with a sufficiently entangled
periodic $\Gamma$ are easily found. However, if we choose instead
another period cell $\PP$ with a smooth boundary, which is not a
planar slab and for which the property (c) is valid, the argument
modifies easily and the claim of Theorem~\ref{leakymain} remains
true. \\ [.2em]
(ii) In the case $d=2$ such a ``puzzle-like'' decomposition can be
always found. To see that, fold the plane into a cylinder of
radius $|b|/2\pi$ so $\Gamma$ becomes a loop which encircles the
cylindrical surface, dividing into two parts $\CC^\pm$ which are
disjoint apart of the common boundary; each of them is connected
because $\Gamma$ has by assumption no self-intersections. Choosing
a point at $\Gamma$, one can thus find two smooth semi-infinite
curves in $\CC^\pm$, even straight from some point on, which go
the two cylinder ``infinities'' without crossing $\Gamma$.
\\ [.2em]
(iii) On the other hand, an analogous decomposition into
translates of a suitable $\PP$ satisfying the hypothesis (c) may
not exist if $d=3$. It depends on the topology of $\Gamma$, a
simple counterexample is given by a ``crotchet-shaped'' curve
which enters $\PP$ on its ``left side'' twice and leaves it once,
and vice versa on the right, without being topologically
equivalent to a line. We conjecture that the claim of
Theorem~\ref{leakymain} remains valid in such situations too,
however, a different method is required to demonstrate it. }
\end{remarks}


\subsection*{Acknowledgments}

P.E. is grateful for the hospitality in Centre de Physique
Th\'eorique, C.N.R.S., where this work was started, and all the
authors to Institut ``Simion Stoilow'' of the Romanian Academy,
where it was finished. The research has been partially supported
by GAAS under the contract A1048101.

\end{document}